\documentclass[conference]{IEEEtran}
\newtheorem{theorem}{Theorem}
\newtheorem{lemma}{Lemma}
\newcommand{\e}{\varepsilon}

\usepackage{cite}

\usepackage{amsmath}
\usepackage{amssymb}

\usepackage{algorithm}
\usepackage{algpseudocode}

\begin{document}

\title{Adaptive Stochastic Mirror Descent \\ for Constrained Optimization}

\author{\IEEEauthorblockN{Anastasia Bayandina}
\IEEEauthorblockA{Moscow Institute of Physics and Technology\\
Moscow, Russia\\
Email: anast.bayandina@gmail.com}}

\maketitle

\begin{abstract}
Mirror Descent (MD) is a well-known method of solving non-smooth convex optimization problems. This paper analyzes the stochastic variant of MD with adaptive stepsizes. Its convergence on average is shown to be faster than with the fixed stepsizes and optimal in terms of lower bounds.
\end{abstract}

\IEEEpeerreviewmaketitle

\section{Introduction}

Recently, the great interest in the optimization problems of a huge scale has been aroused \cite{bib_hugescale}. The dimensions of such problems are of order $10^8$ - $10^{12}$, and computational capabilities only allow for the simplest vector operations. Moreover, it is typical for the corresponding matrices to be sparse. The typical applications are problems related to the Internet \cite{bib_pagerank}, telecommunications, truss topology design \cite{bib_Shpirko}, \cite{bib_ttd}. If the subgradients are also sparse, the deterministic methods \cite{bib_Polyak}, \cite{bib_Shor} can benefit from the sparsity and be applied for solving huge-scale problems. However, it is not always true. Here, the randomized methods play an important part.

Mirror Descent (MD) method \cite{bib_Yudin}, \cite{bib_Beck} is one of the well-known methods of the non-smooth convex optimization. Its randomized version with adaptive stepsizes is studied in this paper. The problem of interest is minimizing a convex function with convex functional constraints. In the assumption of the stochastic first-order oracle, we prove MD to have optimal convergence rate \cite{bib_Yudin}. Due to the adaptive choice of the steplengths, the rate is no more dependent on the global Lipschitz constant but rather on the average of those obtained at the current points.

The paper is divided into five sections. In Section II we state the problem and introduce the notation; the proposed algorithm is described in Section III with the proof of its rate of convergence on average; Section IV contains the remarks and illustrates the discussed approach with the example.


\section{Preliminaries and Problem Statement}
Let $E$ be the $n$-dimensional vector space. Let $\lVert\cdot\rVert$ be an arbitrary norm in $E$ and $\lVert\cdot\rVert_*$ be the conjugate norm in $E^*$:
$$\lVert\xi\rVert_* = \max\limits_{x} \big\{ \langle \xi, x \rangle, \lVert x \rVert \leq 1 \big\}.$$

Let $\mathcal{X} \subset E$ be a closed convex set. We consider the two convex functions $f : \mathcal{X} \rightarrow \mathbb{R}$ and $g : \mathcal{X} \rightarrow \mathbb{R}$ to be subdifferentiable, i.e.
\begin{gather*}
\forall x, y \in \mathcal{X} \hspace{0.2cm} \exists \nabla f(x) : f(y) \geq f(x) + \langle \nabla f(x), y-x \rangle, \\
\forall x, y \in \mathcal{X} \hspace{0.2cm} \exists \nabla g(x) : g(y) \geq g(x) + \langle \nabla g(x), y-x \rangle.
\end{gather*}

We focus on the problem expressed in the form
\begin{equation}
\label{problem_statement_f}
    f(x) \rightarrow \min\limits_{x\in \mathcal{X}},
\end{equation}
\begin{equation}
\label{problem_statement_g}
    \mathrm{s.t.} \hspace{0.3cm} g(x) \leq 0.
\end{equation}
Denote $x_*$ to be the genuine solution of the problem~(\ref{problem_statement_f}),~(\ref{problem_statement_g}).

Assume that we are equipped with the first-order stochastic oracle \cite{bib_Shapiro}, which given the point $x^k\in \mathcal{X}$ returns the values of $\nabla f(x^k, \xi^k), \nabla g(x^k, \xi^k),$ and $g(x^k)$, where $\xi^1, \xi^2, \dots$ are i.i.d random variables that represent oracle noises. The oracle is unbiased, i.e.
\begin{equation}
    \begin{gathered}
    \mathbb{E}\big[ \nabla f(x^k, \xi^k) \big| \xi^{k-1}, \dots, \xi^1 \big]  = \nabla f(x^k) \in \partial f(x^k), \\
    \mathbb{E}\big[ \nabla g(x^k, \xi^k) \big| \xi^{k-1}, \dots, \xi^1 \big] = \nabla g(x^k) \in \partial g(x^k),
    \end{gathered}\label{grad}
\end{equation}
where $x^k = x^k (\xi^{k-1}, \dots, \xi^1)$, and
$$\lVert \nabla f(x^k, \xi^k) \rVert_* < \infty, \hspace{0.3cm} \lVert \nabla g(x^k, \xi^k) \rVert_* < \infty.$$

Consider $d : \mathcal{X} \rightarrow \mathbb{R}$ to be a distance generating function (d.g.f) which is continuously differentiable and strongly convex, modulus 1, w.r.t. the norm $\lVert\cdot\rVert$, i.e.
$$\forall x, y, \in \mathcal{X} \hspace{0.2cm} \langle d'(x) - d'(y), x-y \rangle \geq \lVert x-y \rVert^2.$$
For all $x, y\in \mathcal{X}$ consider the corresponding Bregman divergence
$$V(x, y) = d(y) - d(x) - \langle d'(x), y-x \rangle.$$

Suppose we are given a constant $R$ such that
\begin{equation}\label{R}
    \sup\limits_{x, y\in \mathcal{X}} V(x, y) \leq R^2.
\end{equation}

For all $x\in \mathcal{X}$, $y\in E^*$ define the proximal mapping operator
$$\mathrm{Mirr}_x (y) = \arg\min\limits_{u\in \mathcal{X}} \big\{ \langle y, u \rangle + V(x, u) \big\}.$$
We make the simplicity assumption, which means that $\mathrm{Mirr}_x (y)$ is easily computable.


\section{Stochastic Mirror Descent}

The following algorithm is proposed to solve the problem~(\ref{problem_statement_f}),~(\ref{problem_statement_g}) in the case of stochastic oracle.

\begin{algorithm}[H]
\caption{Adaptive Stochastic Mirror Descent}\label{MD}
\begin{algorithmic}[1]
    \Procedure{MD}{$\e, R, \mathcal{X}$}
        \State $x^1 \gets \arg\min\limits_{x\in \mathcal{X}} d(x)$
        \State initialize the empty set $I$
        \State $k \gets 0$
        \Repeat
            \State $k \gets k + 1$
            \If{$g(x^k) \leq \e$}
                \State $M_k \gets \lVert \nabla f(x^k, \xi^k) \rVert_{*}$
                \State $h_k \gets R \Big(\sum\limits_{i=1}^k M_i^2\Big)^{-1/2}$
                \State $x^{k+1} \gets \mathrm{Mirr}_{x^k}(h_k \nabla f(x^k, \xi^k))$
                \State add $k$ to $I$
            \Else
                \State $M_k \gets \lVert \nabla g(x^k, \xi^k) \rVert_{*}$
                \State $h_k \gets R \Big(\sum\limits_{i=1}^k M_i^2\Big)^{-1/2}$
                \State $x^{k+1} \gets \mathrm{Mirr}_{x^k}(h_k \nabla g(x^k, \xi^k))$
            \EndIf 
        \Until{$\frac{2R}{k} \Big(\sum\limits_{i=1}^{k} M_i^2\Big)^{1/2} \leq \e$}
        \State $N_I \gets |I|$
        \State $\bar{x}^N \gets \frac{1}{N_I} \sum\limits_{k\in I} x^k$
        \State\Return $\bar{x}^N$
    \EndProcedure
\end{algorithmic}
\end{algorithm}

Denote $J = \{ k\in \overline{1, N} \} / I$.

To prove the convergence estimates of Algorithm \ref{MD} we will need the following lemma.

\begin{lemma}
    Let $\alpha_1, \dots, \alpha_N$ be some non-negative sequence. Then
    \begin{equation}\label{lem_1}
        \sum\limits_{k=1}^N \frac{\alpha_k}{\Big(\sum\limits_{i=1}^k \alpha_i\Big)^{1/2}} \leq 2 \Big(\sum\limits_{k=1}^N \alpha_k\Big)^{1/2}.
    \end{equation}
\end{lemma}
\begin{IEEEproof}
    We prove the statement by induction. For $N=1$ it is obvious. Suppose it holds true for $N = n-1$. Denote $S_n = \sum\limits_{k=1}^n \alpha_k$. Then by induction hypothesis
    \begin{gather*}
        \sum\limits_{k=1}^n \frac{\alpha_k}{\sqrt{S_k}} \leq 2\sqrt{S_{n-1}} + \frac{\alpha_n}{\sqrt{S_n}} = \\
        2\sqrt{S_n - \alpha_n} + \frac{\alpha_n}{\sqrt{S_n}} \stackrel{?}{\leq} 2\sqrt{S_n}, \\
        2\sqrt{S_n (S_n - \alpha_n)} \stackrel{?}{\leq} 2S_n - \alpha_n, \\
        2\Big( \frac{S_n}{\alpha} \Big)^2 - 2 \Big( \frac{S_n}{\alpha} \Big) + 1 \stackrel{?}{\geq} 0,
    \end{gather*}
    where the last inequality always holds true.
\end{IEEEproof}

We will also need the following inequality \cite{bib_Duchi} which is the consequence of how the proximal operator is chosen.

\begin{lemma}
    Let $f$ be some convex subdifferentiable function over the convex set $\mathcal{X}$. Let the sequence $\{x^i\}$ be defined by the update
    $$ x^{k+1} = \mathrm{Mirr}_{x^k}(h_k \nabla f(x^k, \xi^k)). $$
    Then, for any $x_* \in \mathcal{X}$
    \begin{equation}
        \begin{gathered}
        f(x^k) - f(x_*) \leq \frac{1}{h_k} \big( V(x^k, x_*) - V(x_{k+1}, x_*) \big) + \\
        \frac{h_k}{2} \lVert \nabla f(x^k, \xi^k) \rVert^2_* + \langle \nabla f(x^k, \xi^k) - \nabla f(x^k), x^k - x_* \rangle.
        \end{gathered}\label{Duchi}
    \end{equation}
\end{lemma}

The following theorem estimates the convergence rate of Algorithm \ref{MD}.

\begin{theorem}\label{th_1}
    The point $\bar{x}^N$ supplied by Algorithm \ref{MD} satisfies
    \begin{equation}
    \label{xxx}
        \mathbb{E}\big[f(\bar{x}^N)\big] - f(x_{*}) \leq \e, \hspace{0.3cm} g(\bar{x}^N) \leq \e
    \end{equation}
    with the average number of oracle calls equal to
    \begin{equation}\label{exp_N}
        \mathbb{E} \big[ N \big] = \mathbb{E} \Bigg[ \Big\lceil\frac{4M^2 R^2}{\e^2}\Big\rceil \Bigg],
    \end{equation}
    where
    \begin{equation}
        M = \Big( \frac{1}{N} \sum\limits_{k=1}^N M_k^2 \Big)^{1/2}.
    \end{equation}
\end{theorem}
\begin{IEEEproof}
    By the definition of $\bar{x}^N$ and the convexity of $f$,
    \begin{equation}\label{th1_eq1}
        N_I \big( f(\bar{x}^N) - f(x_*) \big) \leq \sum\limits_{k\in I} \big( f(x^k) - f(x_*) \big).
    \end{equation}
    
    Denote
    \begin{equation}\label{delta}
        \delta_k = \begin{cases}
            \langle \nabla f(x^k, \xi^k) - \nabla f(x^k), x^k - x_* \rangle, \text{ if } k\in I, \\
            \langle \nabla g(x^k, \xi^k) - \nabla g(x^k), x^k - x_* \rangle, \text{ if } k\in J.
        \end{cases}
    \end{equation}
    
    Using (\ref{Duchi}), consider the summation
    \begin{gather*}
        \sum\limits_{k\in I} \big( f(x^k) - f(x_*) \big) + \sum\limits_{k\in J} \big( g(x^k) - g(x_*) \big) \leq \\
        \sum\limits_{k=1}^N \frac{1}{h_k} \big( V(x^k, x_*) - V(x^{k+1}, x_*) \big) + \\
        \sum\limits_{k=1}^N \frac{h_k M_k^2}{2} + \sum\limits_{k=1}^N \delta_k.
    \end{gather*}
    As long as by (\ref{R})
    \begin{gather*}
        \sum\limits_{k=1}^N \frac{1}{h_k} \big( V(x^k, x_*) - V(x^{k+1}, x_*) \big) = \frac{1}{h_1} V(x^1, x_*) + \\
        \sum\limits_{k=1}^{N-1} \Big( \frac{1}{h_{k+1}} - \frac{1}{h_k} \Big) V(x^{k+1}, x_*) - \frac{1}{h_N} V(x^N, x_*) \leq \\
        \frac{R^2}{h_1} + R^2 \sum\limits_{k=1}^{N-1} \Big( \frac{1}{h_{k+1}} - \frac{1}{h_k} \Big) = \frac{R^2}{h_N},
    \end{gather*}
    then by the definition of stepsizes and (\ref{lem_1}),
    \begin{gather*}
        \sum\limits_{k\in I} \big( f(x^k) - f(x_*) \big) + \sum\limits_{k\in J} \big( g(x^k) - g(x_*) \big) \leq \\
        \frac{R^2}{h_N} + \sum\limits_{k=1}^N \frac{h_k M_k^2}{2} + \sum\limits_{k=1}^N \delta_k = \\
        R \Big( \sum\limits_{k=1}^N M_k^2 \Big)^{1/2} + \frac{R}{2} \sum\limits_{k=1}^N \frac{M_k^2}{\Big(\sum\limits_{i=1}^k M_i^2\Big)^{1/2}} + \sum\limits_{k=1}^N \delta_k \leq \\
        2R \Big( \sum\limits_{k=1}^N M_k^2 \Big)^{1/2} + \sum\limits_{k=1}^N \delta_k.
    \end{gather*}
    Since for $k\in J$
    $$g(x^k) - g(x_{*}) \geq g(x^k) > \e,$$
    recalling (\ref{th1_eq1}) and the stopping criterion, we get
    \begin{gather*}N_I \big( f(\bar{x}^N) - f(x_*) \big) < \\
    \e N_I - \e N + 2R \Big( \sum\limits_{k=1}^N M_k^2 \Big)^{1/2} + \sum\limits_{k=1}^N \delta_k \leq \e N_I + \sum\limits_{k=1}^N \delta_k.
    \end{gather*}
    As long as the inequality is strict, the case of $N_I = 0$ is impossible. Taking the expectation we obtain
    \begin{equation}
        \mathbb{E}\big[f(\bar{x}^N)\big] - f(x_*) \leq \e + \sum\limits_{k=1}^N \mathbb{E} \Big[ \frac{\delta_k}{N_I} \Big].
    \end{equation}
    
    Let us look at the second term on the right-hand side
    \begin{equation}\label{th1_eq2}
        \mathbb{E} \Big[ \frac{\delta_k}{N_I} \Big] = \mathbb{E} \Big[ \frac{1}{N_I} \mathbb{E} \big[ \delta_k \big| N_I \big]\Big].
    \end{equation}
    We are going to show that $\mathbb{E} \big[ \delta_k \big| N_I \big] = 0$.
    
    First,
    \begin{equation}\label{th1_eq3}
        \mathbb{E} \big[ \delta_k \big| N_I \big] = \mathbb{E} \big[ \delta_k \big| N_I=n \big] \Big|_{n=N_I}.
    \end{equation}
    The last denotation describes how the conditional expectation is taken: we fix $N_I$, then take the expectation of $\delta_k$ and then substitute the fixed value back by the random variable.
    
    Suppose we observe the realization $N_I = n$. Denote by $\mathcal{I}_n$ the set of indices such that
    \begin{equation}\label{I_n}
        k\in \mathcal{I}_n \hspace{0.1cm} \mapsto \hspace{0.1cm} g(x^k) \leq \e.
    \end{equation}
    There are $\binom{N}{n}$ distinct sets $\mathcal{I}_n$ and for each of them $|\mathcal{I}_n| = n$. By the law of total expectation,
    \begin{equation}\label{th1_eq4}
        \mathbb{E} \big[ \delta_k \big| N_I=n \big] \Big|_{n=N_I} = \Big( \sum\limits_{\mathcal{I}_n} \mathbb{E} \big[ \delta_k \big| \mathcal{I}_n \big] \mathbb{P}\big( \mathcal{I}_n \big) \Big) \Big|_{n=N_I}.
    \end{equation}
    Let us show that for each set of indices $\mathcal{I}_n$ defined in (\ref{I_n}) holds $\mathbb{E} \big[ \delta_k \big| \mathcal{I}_n \big] = 0$. Applying the "tower property" of conditional expectation and recalling (\ref{delta}) and (\ref{grad}), we derive
    \begin{equation}
        \mathbb{E} \big[ \delta_k \big| \mathcal{I}_n \big] = \mathbb{E} \Big[ \mathbb{E}\big[ \delta_k \big| \mathcal{I}_n \big] \Big| x^k \Big] = \mathbb{E} \big[ \delta_k \big| x^k \big] = 0.
    \end{equation}
    Thus, we have proved that together with (\ref{th1_eq4}), (\ref{th1_eq3}), (\ref{th1_eq2})
    \begin{equation}
        \sum\limits_{k=1}^N \mathbb{E} \Big[ \frac{\delta_k}{N_I} \Big] = 0
    \end{equation}
    and the result holds.
    
    For $i\in I$ holds $g(x^i) \leq \e$. Then, by the definition of $\bar{x}^N$ and the convexity of $g$,
    \begin{gather*}
        N_I g(\bar{x}^N) \leq \sum\limits_{k\in I} g(x^k) \leq \e N_I.
    \end{gather*}
\end{IEEEproof}


\section{Remarks}

\textbf{Remark 1}. In the assumption of the uniform boundedness of subgradients, i.e. for each $x\in \mathcal{X}$ and $\xi \in \{\xi^1, \xi^2, \dots\}$
$$\lVert \nabla f(x, \xi) \rVert_*^2 \leq M^2, \hspace{0.1cm} \lVert \nabla g(x, \xi) \rVert_*^2 \leq M^2$$
for some constant $M$, we would have the exact estimate of the number of oracle calls
\begin{equation}\label{N}
    N = \Big\lceil\frac{4M^2 R^2}{\e^2}\Big\rceil.
\end{equation}
Note that (\ref{exp_N}) is no worse than (\ref{N}) and is even better in practice since $M_k$ tend to get smaller as the method approaches to the minimizer.

\textbf{Remark 2}. With the fixed stepsizes the estimation is \cite{bib_Juditsky}
$$N = \Big\lceil\frac{2M^2 R_0^2}{\e^2}\Big\rceil.$$
Here $R_0^2 = V(x^1, x_*)$ which is smaller then $R^2$ in general.

\textbf{Remark 3}. Common applications seek to solve problems in the form
\begin{equation}
    f(x) = \frac{1}{2} \langle x, Ax \rangle \rightarrow \min\limits_{x\in S_n(1)},
\end{equation}
\begin{equation}\label{example_g}
    \mathrm{s.t.} \hspace{0.3cm} g(x) = \max\limits_{m\in \overline{1, \mathcal{M}}} \big\{ \langle c_m, x \rangle \big\} \leq 0,
\end{equation}
where $S_n(1)$ is the $n$-dimensional simplex, $A$ is some $n\times n$ matrix and $\big\{ c_m \big\}_{m=1}^{\mathcal{M}}$ are some vectors in~$\mathbb{R}^n$.
        
Firstly, it is worth mentioning that typically the vectors $\big\{ c_m \big\}_{m=1}^{\mathcal{M}}$ are sparse. For the constraints of the type (\ref{example_g}) it means that their subgradients are also sparse. Then, each update in MD changes just a few entries in $x^k$. Choosing the appropriate distribution of oracle noise we can reduce the update to changing only one entry in the vector \cite{bib_hugescale}.

Secondly, even if the matrix $A$ is sparse, the gradient $\nabla f(x) = Ax$ is usually not. The exact computation of the gradient takes $O(n^2)$ arithmetic operations, which is bad for huge-scale optimization problems. However, the possible way out of this situation can be randomization \cite{bib_JuNem}. We generate the appropriate distribution over the columns of the matrix $A$ and take just one column, according to the distribution, to be the unbiased estimate of the gradient. Let us construct the distribution and show that the estimate is unbiased.

Let $\xi$ be a random variable taking its values in $\{1, \dots, n\}$ and let $A^{\langle i \rangle}$ denote the $i$th column of the matrix $A$. Suppose the current point is $x = (x_1, \dots, x_n)$. Since $x\in S_n(1)$,
\begin{gather*}\mathbb{E}\big[ A^{\langle \xi \rangle} \big] = A^{\langle 1 \rangle} \underbrace{\mathbb{P}\big( \xi = 1 \big)}_{x_1} + \dots + A^{\langle n \rangle} \underbrace{\mathbb{P}\big( \xi = n \big)}_{x_n} = \\
A^{\langle 1 \rangle} x_1 + \dots + A^{\langle n \rangle} x_n = Ax.
\end{gather*}

In this case, the computation of a subgradient requires only $O(n)$ arithmetic operations.


\section{Conclusion}
We studied the stochastic version of MD for minimizing convex functions with convex functional constraints. It has been proposed to modify MD so that the steplengths become adaptive and depend only on the subgradient values at current points. Furthermore, it has been proved to converge practically faster than the method with fixed steplengths because of the improved constant, which also has become adaptive. Eventually, the method has been briefly illustrated with a possible application.


\section*{Acknowledgment}

The author gratefully acknowledges the help and valuable discussion kindly provided by Dr. Gasnikov.

This research was funded by Russian Science Foundation (project 17-11-01027).


\end{document}